\documentclass[12pt]{article}

\newcommand{\open}[1]{\left#1 \vphantom{ \vert^2 } \right.}
\newcommand{\close}[1]{\left. \vphantom{ \vert^2 } \right#1}
\renewcommand{\span}{{\rm span}}
\newcommand{\inverse}[1]{{\textstyle\frac{1}{#1}}}
\newcommand{\half}{\inverse{2}}
\newcommand{\quarter}{\inverse{4}}
\newcommand{\uq}{\underline{q}}
\newcommand{\qed}{\hfill$\Box$}

\title{{\Large \bf The direct $L^2$ geometric structure on a manifold of probability densities with applications to Filtering}}
\author{Damiano Brigo\thanks{I am grateful towards Giuseppe Tinaglia and Alexander Pushnitski
 for help with geometry and topology.  All remaining errors are my own. }  \\
Dept. of Mathematics\\
King's College, London\\
{\tt damiano.brigo@kcl.ac.uk} }

\date{\small {First Version Posted on Nov 12, 2011, on hal.archives-ouvertes.fr and damianobrigo.it \\ Second Version Posted on Nov 27, 2011, on arXiv.org. This version: \today }}
\newtheorem{theorem}{Theorem}[section]

\newtheorem{definition}[theorem]{Definition}
\newtheorem{example}[theorem]{Example}
\newtheorem{remark}[theorem]{Remark}
\newtheorem{problem}[theorem]{Problem}
\begin{document}
%

\maketitle
\thispagestyle{empty}
\begin{abstract}
In this paper we introduce a projection method for the space of probability distributions based on the differential geometric approach to statistics. This method is based on a direct $L^2$ metric as opposed to the usual Hellinger distance and the related Fisher Information metric. We explain how this apparatus can be used for the nonlinear filtering problem, in relationship also to earlier projection methods based on the Fisher metric. Past projection filters focused on the Fisher metric and the exponential families that made the filter correction step exact. In this work we introduce the mixture projection filter, namely the projection filter based on the direct $L^2$ metric and based on a manifold given by a mixture of pre-assigned densities. The resulting prediction step in the filtering problem is described by a linear differential equation, while the correction step can be made exact. We analyze the relationship of a specific class of $L^2$ filters with the Galerkin based nonlinear filters, and highlight the differences with our approach, concerning particularly the continuous--time observations filtering problems.
\end{abstract}

\bigskip


\bigskip

\noindent {\bf AMS Classification codes: 53B25, 53B50, 60G35, 62E17, 62M20, 93E11}

\pagestyle{myheadings}
\markboth{}
{D. Brigo. The direct $L^2$ geometric structure with applications to Filtering}

\section{Introduction}
In this paper we consider the nonlinear filtering problem
in continuous time.
For a quick introduction to the filtering problem see
Davis and Marcus (1981) \cite{davis81b}. For a more complete
treatment see Liptser and Shiryayev (1978) \cite{LipShi}
from a mathematical point of view or Jazwinski (1970) \cite{jazwinski70a}
for a more applied perspective. For recent results see the collection of papers \cite{crisan}.

The nonlinear filtering problem has an infinite--dimensional solution
in general. Constructing of approximate finite-dimensional filters is an important area of research.

When the system has continuous time signal and continuous time observations, the solution of the filtering problem is a Stochastic PDE which can be seen as a generalization of the Fokker--Planck equation expressing the evolution of the density of a diffusion process. This filtering equation is called Kushner--Stratonovich equation, and an unnormalized (simpler) version of it is known as the Duncan--Mortensen--Zakai Stochastic Partial Differential
Equation. When observations are in discrete time, the filtering problem decomposes into a prediction step, given by the Fokker-Planck equation, and a correction step, given by Bayes formula.

In \cite{BrigoPhD}, \cite{brigo98} and \cite{brigo99} the Fisher metric is used to project the Kushner--Stratonovich (or the Fokker--Planck) equation onto an exponential family of probability densities, yielding the new class of approximate filters called {\em projection filters}. The projection filters are based on the differential geometric approach to statistics, as developed by \cite{amari85a} and \cite{pistonesempi}.
It is also shown that one can choose the family so as to make the prediction step exact. Moreover, it is shown that for exponential families the projection filters coincide with the assumed density filters.

In \cite{Brigo2,Brigo3} the Gaussian projection filter is studied in the small-noise setting.

In the present paper we choose a different differential geometric structure based on a direct $L^2$ metric as opposed to the usual Hellinger distance and the related Fisher Information metric. We explain how this structure can be used to derive a different family of finite dimensional filters that form a good approximation for the solution of the nonlinear filtering problem. This structure is particularly suited to be applied to mixture families of distributions, similarly to how exponential families are well suited to work with the Fisher information metric.  In this work we thus introduce the mixture projection filter, namely the projection filter based on the direct $L^2$ metric and based on a manifold given by a mixture of pre-assigned densities. One key result we obtain is that the prediction step is given by a linear differential equation, whereas the correction step can be made exact by updating the basis functions for the tangent space of the manifold, namely the mixture components, at each observation time.

The exponential projection filter had a clear relationship with the assumed density filters, as documented in \cite{brigo99}. The $L^2$ mixture projection filter has a clear relationship with earlier Galerkin-based approaches to non-linear filtering, see for example \cite{beard} and \cite{kenney}. However, the geometric structure and the exact $L^2$ projection make the method in this paper more general, giving the possibility to apply it to manifolds that are more general than the standard mixture family. Morevoer, in the continuous time observations case, the $L^2$ projection filter is based on the Stratonovich calculus that is needed to keep the projected dynamics into the tangent space of the manifold, whereas the Galerkin based projection filter \cite{beard} is based on Ito calculus. We will explore in detail our mixture projection filter based on direct $L^2$ metric as compared with the Galerkin methods in future research, where we will also implement the mixture projection filter equations numerically, both for the simple mixture family and for more general families to which the Galerkin method cannot be applied. We will also investigate the choice of the specific mixture family, starting with gaussian or lognormal mixtures, with possible likelihood ratio corrections that make the correction step exact and allow for the definition of a rigorous measure for the filtering error, based on a projection residual.

\section{Statistical manifolds} \label{StatMan}

On the measurable space $({\bf R}^n,{\cal B}({\bf R}^n))$
we consider a non--negative and $\sigma$--finite measure $\lambda$,
and we define ${\cal M}(\lambda)$ to be the set of all non--negative
and finite measures $\mu$ which are absolutely continuous
w.r.t.\ $\lambda$, and whose density
\begin{displaymath}
   p_{\mu} = \frac{d \mu}{d \lambda}
\end{displaymath}
is positive $\lambda$--a.e.
For simplicity, we restrict ourselves to the case
where $\lambda$ is the Lebesgue measure on ${\bf R}^n$.
We also assume that the total measure is normalized to one, so as to represent a probability measure. This in turn implies that $p_\mu$ integrates to one.

In the following, we denote by $H(\lambda)$ the set of all the densities
of measures contained in  ${\cal M}(\lambda)$.
Notice that, as all the measures in ${\cal M}(\lambda)$ are
non--negative and finite,
we have that if $p$ is a density in $H(\lambda)$
then $p\in L_1(\lambda)$, that is $\sqrt{p}\in L^2(\lambda)$.
The above remark implies that
the set ${\cal R}(\lambda):=\{\sqrt{p}:p{\in}H(\lambda)\}$
of square roots of densities of $H(\lambda)$ is a subset of $L^2(\lambda)$.
Notice that all $\sqrt{p}$ in ${\cal R}(\lambda)$
satisfy $\sqrt{p(x)}>0$, for almost every $x\in {\bf R}^n$.

We notice the important point that neither $H(\lambda)$ nor ${\cal R}(\lambda)$
are vector subspaces of $L_1$ or $L^2$ respectively. Hence, we cannot view them as normed subspaces or topological vector spaces.

We will be able to use the $L^2$ norm to define a {\em metric} in ${\cal R}$, but we will not be able to view ${\cal R}$ as a normed space.

\subsection{The Hellinger distance}

The above remarks lead to the definition of the following metric
in ${\cal R}(\lambda)$, see Jacod and Shiryayev~\cite{jacod87a}
or Hanzon~\cite{hanzon89a},
$d_{\cal R}(\sqrt{p_1},\sqrt{p_2}) := \Vert \sqrt{p_1}-\sqrt{p_2} \Vert$,where $\Vert\cdot\Vert$ denotes the norm of the Hilbert space $L^2(\lambda)$.
This leads to the Hellinger metric on $H(\lambda)$
(or ${\cal M}(\lambda)$), obtained
by using the bijection between densities (or measures) and
square roots of densities~: if $\mu_1$ and $\mu_2$ are the measures
having densities $p_1$ and $p_2$ w.r.t.\ $\lambda$,
the Hellinger metric is defined
as $d_{\cal M}(\mu_1,\mu_2)=d_H(p_1,p_2)=d_{\cal R}(\sqrt{p_1},\sqrt{p_2})$.
It can be shown, see e.g.~\cite{hanzon89a}, that the
distance $d_{\cal M}(\mu_1,\mu_2)$ in ${\cal M}(\lambda)$ is defined
independently of the particular $\lambda$ we choose as basic measure,
as long as both $\mu_1$ and $\mu_2$ are absolutely continuous
w.r.t.\ $\lambda$.
As one can always find a $\lambda$ such that both $\mu_1$
and $\mu_2$ are absolutely continuous w.r.t.\ $\lambda$
(take for example $\lambda:=(\mu_1+\mu_2)/2$),
the distance is well defined on the set of all finite
and positive measures on $(\Omega,{\cal F})$.

\subsection{The $L^2$ direct distance}
There is another possibility for defining a metric in $H$.
We consider the following subset of $H$:
\[ H_2(\lambda) = H(\lambda) \cap L^2(\lambda) \]
i.e. the set of $L^2$ densities. Notice that here we do not take the square root, but we use the $L^2$ structure directly on the densities. If we further assume that densities in $H$ are bounded, then
\[ H_2(\lambda) = H(\lambda)  \]
since bounded positive functions that are in $L_1$ are also in $L^2$.

This structure leads to the definition of the following metric
in $H_2(\lambda)$:
$d_{2}(p_1,p_2) := \Vert p_1-p_2 \Vert$. $H_2$ with this metric is a metric space but, again, it is not a normed space, since it is not a vector space. We call this metric the direct $L^2$ distance, since it is taken directly on the densities rather than mapping them to their square roots.


\subsection{Neither $(H(\lambda),d_H)$ nor $(H_2(\lambda),d_2)$ are $L^2$ Hilbert manifolds}

Despite being subsets of $L^2$, neither $(H(\lambda),d_H)$ (or the equivalent $({\cal R}(\lambda),d_{\cal R})$) nor $(H_2(\lambda),d_2)$ are locally homeomorphic to $L^2(\lambda)$, hence they are not manifolds modeled on $L^2(\lambda)$.
Indeed, any open set of $L^2(\lambda)$ contains functions which are negative in a set with positive $\lambda$--measure.
There is no open set of $L^2(\lambda)$ which contains only positive functions such as the functions of $H_2(\lambda)$ or ${\cal R}(\lambda)$.

\subsection{Definition of Tangent vectors through the $L^2$ structure}

Consider an open subset $M$ of $L^2(\lambda)$.
Let $x$ be a point of $M$, and let
$\gamma:(-\epsilon,\epsilon)~\rightarrow~M$ be a curve on $M$
around $x$, i.e.\ a differentiable map
between an open neighborhood of $0\in {\bf R}$ and $M$ such
that $\gamma(0) = x$. We can define the tangent vector to $\gamma$
at $x$ as the Fr\'echet derivative
$D\gamma(0):(-\epsilon,\epsilon)~\rightarrow~L^2(\lambda)$,
i.e.\ the linear map
defined in ${\bf R}$ around $0$ and taking values in $L^2(\lambda)$
such that the following limit holds~:
\begin{displaymath}
   \lim_{\vert h \vert \rightarrow 0}\;
   \frac{\Vert \gamma(h)-\gamma(0)-D\gamma(0)\cdot h \Vert}{
   \vert h \vert} = 0\ .
\end{displaymath}
The map $D\gamma(0)$ approximates linearly the change of $\gamma$
around $x$.
Let ${\cal C}_x(M)$ be the set of all the curves on $M$ around $x$.
If we consider the space
\begin{displaymath}
   L_x M := \{ D\gamma(0)\,:\,\gamma\in {\cal C}_x(M) \}\ ,
\end{displaymath}
of tangent vectors to all the possible curves on $M$ around $x$,
we obtain again the space $L^2(\lambda)$. This is due to the
fact that for every $v \in L^2(\lambda)$ we can always consider the
straight line $\gamma^v(h) := x + h\, v$. Since $M$ is open, $\gamma^v(h)$
takes values in $M$ for $\vert h \vert$ small enough.
Of course $D\gamma^v(0) = v$, so that indeed $L_x M = L^2(\lambda)$.

\subsection{Finite dimensional submanifold embedded in $L^2$}

The situation becomes different if we consider an {\em $m$--dimensional}
manifold $N$ that is a subset of $L^2$ (and, possibly, a subset of ${\cal R}$ or $H_2$ above).
As such, it can be endowed with the topology induced by the $L^2$ norm.
Because $N$ is $m$-dimensional, it is also locally homeomorphic to $R^m$.


We can consider the induced $L^2$ structure on $N$ as follows~:
suppose $x \in N$, and define again
\begin{displaymath}
   L_x N := \{ D\gamma(0)\,:\,\gamma\in {\cal C}_x(N) \}\ .
\end{displaymath}
This is a linear subspace of $L^2(\lambda)$
called the {\em tangent vector space} at $x$, which does not coincide
with $L^2(\lambda)$ in general (due to the finite dimension of $N$, this tangent space will be $m$-dimensional).
The set of all tangent vectors at all points $x$ of $N$
is called the {\em tangent bundle}, and will be denoted by $L\, N$.
In our work we shall consider finite dimensional manifolds $N$
embedded in $L^2(\lambda)$,
which are contained in ${\cal R}(\lambda)$ or $H_2$ as a set,
i.e.\ $N \subset {\cal R}(\lambda) \subset L^2(\lambda)$ or $N \subset H_2(\lambda) \subset L^2(\lambda)$,
so that usually $x=\sqrt{p}$ or $x=p$, respectively.

We analyze the two cases separately.

\subsection{Finite dimensional manifolds $N$ in $({\cal R},d_{\cal R})$}

If $N$ is $m$--dimensional, it is locally homeomorphic to ${\bf R}^m$,
and it may be described locally by a chart~: if $\sqrt{p}\in N$,
there exists a pair $(S^{1/2},\phi)$ with $S^{1/2}$ open neighbourhood
of $\sqrt{p}$ in $N$ for the topology induced by $d_{\cal R}$ and $\phi:S^{1/2}\rightarrow\Theta$ homeomorphism
of $S^{1/2}$ with the topology induced by $d_{\cal R}$ onto an open subset $\Theta$ of ${\bf R}^m$ with the usual topology of ${\bf R}^m$.
%
%
%
%
By considering the inverse map $i$ of $\phi$,
\begin{eqnarray*}
   i : \Theta &\longrightarrow& S^{1/2} \\
   \theta &\longmapsto& \sqrt{p(\cdot,\theta)}
\end{eqnarray*}
we can express $S^{1/2}$ as
\begin{displaymath}
   i(\Theta) = \{ \sqrt{p(\cdot,\theta)}\,,\,\theta\in \Theta \}
   = S^{1/2}.
\end{displaymath}

We will work only with the single coordinate chart $(S^{1/2},\phi)$
as it is done in~\cite{amari85a}.
From the fact that $(S^{1/2},\phi)$ is a chart, it follows that
\begin{displaymath}
   \{ \frac{\partial i(\cdot,\theta)}{\partial \theta_1},\cdots,
   \frac{\partial i(\cdot,\theta)}{\partial \theta_m} \}
\end{displaymath}
is a set of linearly independent vectors in $L^2(\lambda)$.
In such a context, let us see what
the vectors of $L_{\sqrt{p(\cdot,\theta)}}S^{1/2}$ are.
We can consider a curve in $S^{1/2}$ around $\sqrt{p(\cdot,\theta)}$
to be of the form $\gamma:h{\mapsto}\sqrt{p(\cdot,\theta(h))}$,
where $h\mapsto\theta(h)$ is a curve in $\Theta$ around $\theta$.
Then, according to the chain rule,
we compute the following Fr\'echet derivative:
\begin{displaymath}
   D \gamma(0) = \left. D \sqrt{p(\cdot,\theta(h))} \right\vert_{h=0}
   = \sum_{k=1}^m \frac{\partial \sqrt{p(\cdot,\theta)}}{
   \partial \theta_k}\; \dot{\theta}_k(0)
   = \sum_{k=1}^m \frac{1}{2 \sqrt{p(\cdot,\theta)}}
   \frac{\partial p(\cdot,\theta)}{\partial \theta_k}\;
   \dot{\theta}_k(0)\ .
\end{displaymath}
We obtain that a basis for the tangent vector space
at $\sqrt{p(\cdot,\theta)}$
to the space $S^{1/2}$ of square roots of densities of $S$
is given by~:
\begin{equation} \label{sqrtbas}
   L_{\sqrt{p(\cdot,\theta)}}S^{1/2}
   = \span \{\frac{1}{2 \sqrt{p(\cdot,\theta)}}
   \frac{\partial p(\cdot,\theta)}{\partial \theta_1},\cdots,
   \frac{1}{2 \sqrt{p(\cdot,\theta)}}
   \frac{\partial p(\cdot,\theta)}{\partial \theta_m} \}\ .
\end{equation}
As $i$ is the inverse of a chart, these vectors are actually linearly
independent, and they indeed form a basis of the tangent vector space.
One has to be careful, because if this were not true, the dimension of
the above spanned space could drop.


The inner product of any two basis elements is defined,
according to the $L^2$ inner product
\begin{eqnarray} \label{hels}
   \langle \frac{1}{2 \sqrt{p(\cdot,\theta)}}\,
   \frac{\partial p(\cdot,\theta)}{\partial \theta_i},
   \frac{1}{2 \sqrt{p(\cdot,\theta)}}\,
   \frac{\partial p(\cdot,\theta)}{\partial \theta_j} \rangle
  =
   \quarter \int \frac{1}{p(x,\theta)}\,
   \frac{\partial p(x,\theta)}{\partial \theta_i}\,
   \frac{\partial p(x,\theta)}{\partial \theta_j}\,
   d\lambda(x)
   = \quarter\; g_{ij}(\theta)\ .
\end{eqnarray}
This is, up to the numeric factor~$\quarter$, the Fisher information metric,
see for example~\cite{amari85a}, \cite{murray93a}
and~\cite{aggrawal74a}.
The matrix $g(\theta) = (g_{ij}(\theta))$ is called
the Fisher information matrix.

Next, we introduce the orthogonal projection between any linear
subspace $V$ of $L^2(\lambda)$
containing the finite dimensional tangent vector space~(\ref{sqrtbas})
and the tangent vector space~(\ref{sqrtbas}) itself.
Let us remember that our basis is not orthogonal, so that we have
to project according to the following formula:
\begin{displaymath}
\begin{array}{rcl}
   \Pi~: V &\longrightarrow& \span \{w_1,\cdots,w_m \} \\ \\
   v &\longmapsto& \displaystyle \sum_{i=1}^m [ \sum_{j=1}^m
   W^{ij}\; \langle v,w_j \rangle ]\; w_i
\end{array}
\end{displaymath}
where $\{w_1,\cdots,w_m\}$ are $m$ linearly independent
vectors, $W:=(\langle w_i,w_j \rangle)$ is the matrix
formed by all the possible inner products of such linearly independent
vectors, and $(W^{ij})$ is the inverse of the matrix $W$.
In our context $\{w_1,\cdots,w_m\}$ are the vectors in (\ref{sqrtbas}),
and of course $W$ is, up to the numeric factor~$\quarter$,
the Fisher information matrix given by~(\ref{hels}).
Then we obtain the following projection formula,
where $(g^{ij}(\theta))$ is the inverse of the Fisher
information matrix $(g_{ij}(\theta))$~:
\begin{equation} \label{proL2}
\begin{array}{rcl}
  && \Pi_\theta~: L^2(\lambda)\supseteq V \longrightarrow \displaystyle
   \span \{ \frac{1}{2 \sqrt{p(\cdot,\theta)}}\,
   \frac{\partial p(\cdot,\theta)}{\partial \theta_1},\cdots,
   \frac{1}{2 \sqrt{p(\cdot,\theta)}}\,
   \frac{\partial p(\cdot,\theta)}{\partial \theta_m} \} \\ \\
  && \Pi_\theta[v] = \displaystyle \sum_{i=1}^m [ \sum_{j=1}^m
   4 g^{ij}(\theta)\; \langle v,\frac{1}{2 \sqrt{p(\cdot,\theta)}}\,
   \frac{\partial p(\cdot,\theta)}{\partial \theta_j} \rangle ]\;
   \frac{1}{2 \sqrt{p(\cdot,\theta)}}\,
   \frac{\partial p(\cdot,\theta)}{\partial \theta_i}\ .
\end{array}
\end{equation}
Let us go back to the definition of tangent vectors for our statistical
manifold.
Amari~\cite{amari85a} uses a different representation of tangent vectors
to $S$ at $p$. Without exploring all the assumptions needed, let us
say that Amari defines an isomorphism between the actual tangent space
and the vector space
\begin{displaymath}
   \span \{
   \frac{\partial \log p(\cdot,\theta)}{\partial \theta_1},\cdots,
   \frac{\partial \log p(\cdot,\theta)}{\partial \theta_m} \}\ .
\end{displaymath}
On this representation of the tangent space, Amari defines a Riemannian
metric given by
\begin{displaymath}
   E_{p(\cdot,\theta)}
   \{ \frac{\partial \log p(\cdot,\theta)}{\partial \theta_i}\;
   \frac{\partial \log p(\cdot,\theta)}{\partial\theta_j} \}\ ,
\end{displaymath}
where $E_p\{\cdot\}$ denotes the expectation w.r.t.\
the probability density $p$.
This is again the Fisher information metric, and indeed this is the
most frequent definition of Fisher metric.
In fact, it is easy to check that
\begin{eqnarray}
   && E_{p(\cdot,\theta)}
   \{ \frac{\partial \log p(\cdot,\theta)}{\partial \theta_i}\;
   \frac{\partial \log p(\cdot,\theta)}{\partial \theta_j} \}
   = \int \frac{\partial \log p(x,\theta)}{\partial \theta_i}\,
   \frac{\partial \log p(x,\theta)}{\partial \theta_j}\,
   p(x,\theta)\,d\lambda(x) \nonumber \\ \label{Fish1} \\ \nonumber
   && \qquad = \int \frac{1}{p(x,\theta)}\,
   \frac{\partial p(x,\theta)}{\partial \theta_i}\,
   \frac{\partial p(x,\theta)}{\partial \theta_j}\,
   d\lambda(x) = g_{ij}(\theta)\ .
\end{eqnarray}
From the above relation and from~(\ref{hels})
it is clear that, up to the numeric factor~$\quarter$,
the Fisher information metric and the Hellinger metric coincide on
the two representations of the tangent space to $S$
at $p(\cdot,\theta)$.

There is another way of measuring how close two densities of $S$ are.
Consider the Kullback--Leibler information between two densities $p$
and $q$ of $H(\lambda)$~:
\begin{displaymath}
   K(p,q) := \int \log \frac{p(x)}{q(x)}\; p(x)\, d\lambda(x)
   = E_{p}\{\log \frac{p}{q} \}\ .
\end{displaymath}
This is not a metric, since it is not symmetric and it does
not satisfy the triangular inequality. When applied to a finite
dimensional manifold such as $S$, both the Kullback--Leibler
information and the Hellinger distance are particular cases
of $\alpha$--divergence, see~\cite{amari85a} for the details.
One can show that the Fisher metric and the Kullback--Leibler
information coincide infinitesimally.
Indeed, consider the two densities $p(\cdot,\theta)$
and $p(\cdot,\theta+d\theta)$ of $S$.
By expanding in Taylor series, we obtain
\begin{eqnarray*}
   K(p(\cdot,\theta),p(\cdot,\theta+d\theta))
   &=& - \sum_{i=1}^m E_{p(\cdot,\theta)}\{
   \frac{\partial \log p(\cdot,\theta)}{\partial \theta_i} \}\,
   d\theta_i \\ \\
   && \qquad - \sum_{i,j=1}^m E_{p(\cdot,\theta)}\{
   \frac{\partial^2 \log p(\cdot,\theta)}{
   \partial \theta_i\partial \theta_j} \}\, d\theta_i\, d\theta_j
   + O(\vert d\theta \vert^3) \\ \\
   &=& \sum_{i,j=1}^m g_{ij}(\theta)\, d\theta_i\, d\theta_j
   + O(\vert d\theta \vert^3)\ .
\end{eqnarray*}
The interested reader is referred to~\cite{aggrawal74a}.

\begin{example} {\bf (The Gaussian family and the Fisher metric with canonical parameters)}. We may consider the Fisher metric for the Gaussian family of densities. The Gaussian family may be defined as a particular exponential family, represented with canonical parameters $\theta$, given by
\[ \{ p(x,\theta) = \exp( \theta_1 x + \theta_2 x^2 - \psi(\theta)), \ \theta_2 <0\} \]
where one has easily
\[ \psi(\theta) = \half \ln\left(\frac{\pi}{-\theta_2}\right) - \frac{\theta_1^2}{4 \theta_2}\]
and the Fisher metric is
\[ g(\theta) = \left[
\begin{array}{cc} -1/(2 \theta_2) & \theta_1 / (2 \theta_2^2)\\  \theta_1 / (2 \theta_2^2) &
1/(2 \theta_2^2) - \theta_1^2/(2 \theta_2^3)
\end{array} \right]
\]
The familiar representation of Gaussian densities is in terms of mean and variance, given respectively by
\[ \mu = - \theta_1/(2 \theta_2),  \ \ v = \sigma^2 = (1/\theta_2 - \theta_1^2/\theta_2^2)/2 \]
\end{example}

The Fisher metric is used ideally to compute the distance between two infinitesimally near points $p(\cdot,\theta)$ and $p(\cdot,\theta+ d \theta)$. Informally, we can write
\[ d_H(p(\cdot,\theta),p(\cdot,\theta+ d \theta)) = (d \theta)^T g(\theta) d \theta   \]
Notice that the matrix changes when changing coordinates, whereas the distance must clearly be the same. Hence if we have another set of coordinates $\eta$ related by diffeomorphism $\eta = \eta(\theta)$ to $\theta$, with inverse $\theta = \theta(\eta)$, then clearly

\[ d_H(p(\cdot,\eta),p(\cdot,\eta+ d \eta)) = (d \eta)^T \  (\partial_\eta \theta(\eta))^T \  g(\theta(\eta)) \  \partial_\eta \theta(\eta)\  d \eta \]
where $\partial_\eta \theta(\eta)$ is the Jacobian matrix of the transformation. It follows that
\[ g(\eta) = (\partial_\eta \theta(\eta))^T \  g(\theta(\eta)) \  \partial_\eta \theta(\eta) \]
\begin{example} {\bf (The Gaussian family and the Fisher metric with expectation parameters)}. We may consider the Fisher metric for the Gaussian family of densities in the parameters $\mu$ and $v$. These are related to the so called expectation parameters $\mu$ and $v + \mu^2$. With this coordinate system the Fisher metric is much simpler and the matrix is diagonal, resulting in
\[  g(\mu,v) = \frac{1}{v} \left[
\begin{array}{cc} 1 & 0 \\  0 &
1/(2 v)
\end{array} \right]
\]
This can be derived either by applying the change of coordinates formula, or Eq. \ref{hels} directly, with the parameters $\theta_1,\theta_2$ replaced by $\mu,v$.
\end{example}

\subsection{Finite dimensional manifolds $N$ in $(H_2,d_2)$}

Alternatively, if we use $H_2$ instead of ${\cal R}$ as a set where $N$ is contained,
$N$ can still be described locally by a chart~: if ${p}\in N$,
there exists a pair $(S,\psi)$ with $S$ open neighbourhood
of ${p}$ in $N$ for the topology induced by $d_2$ and $\psi:S\rightarrow\Theta$ homeomorphism
of $S$ with the topology induced by $d_2$ onto an open subset $\Theta$ of ${\bf R}^m$ with the usual topology.

By considering the inverse map $j$ of $\psi$,
\begin{eqnarray*}
   j : \Theta &\longrightarrow& S \\
   \theta &\longmapsto& {p(\cdot,\theta)}
\end{eqnarray*}
we can express $S$ as
\begin{displaymath}
   j(\Theta) = \{ {p(\cdot,\theta)}\,,\,\theta\in \Theta \}
   = S.
\end{displaymath}

We will work only with the single coordinate chart $(S,\psi)$.
From the fact that $(S,\psi)$ is a chart, it follows that
\begin{displaymath}
   \{ \frac{\partial j(\cdot,\theta)}{\partial \theta_1},\cdots,
   \frac{\partial j(\cdot,\theta)}{\partial \theta_m} \}
\end{displaymath}
is a set of linearly independent vectors in $L^2(\lambda)$.
In such a context, let us see what
the vectors of $L_{{p(\cdot,\theta)}}S$ are.
We can consider a curve in $S$ around ${p(\cdot,\theta)}$
to be of the form $\gamma:h{\mapsto}{p(\cdot,\theta(h))}$,
where $h\mapsto\theta(h)$ is a curve in $\Theta$ around $\theta$.
Then, according to the chain rule,
we compute the following Fr\'echet derivative:
\begin{displaymath}
   D \gamma(0) = \left. D {p(\cdot,\theta(h))} \right\vert_{h=0}
   = \sum_{k=1}^m \frac{\partial {p(\cdot,\theta)}}{
   \partial \theta_k}\; \dot{\theta}_k(0)
   = \sum_{k=1}^m
   \frac{\partial p(\cdot,\theta)}{\partial \theta_k}\;
   \dot{\theta}_k(0)\ .
\end{displaymath}
We obtain that a basis for the tangent vector space
at ${p(\cdot,\theta)}$
to the space $S$
is given by~:
\begin{equation} \label{sqrtbas2}
   L_{{p(\cdot,\theta)}}S
   = \span \{
   \frac{\partial p(\cdot,\theta)}{\partial \theta_1},\cdots,
   \frac{\partial p(\cdot,\theta)}{\partial \theta_m} \}\ .
\end{equation}
As $j$ is the inverse of a chart, these vectors are actually linearly
independent, and they indeed form a basis of the tangent vector space.
One has to be careful, because if this were not true, the dimension of
the above spanned space could drop.


The inner product of any two basis elements is defined,
according to the $L^2$ inner product
\begin{eqnarray} \label{SImetric}
   \langle
   \frac{\partial p(\cdot,\theta)}{\partial \theta_i},
   \frac{\partial p(\cdot,\theta)}{\partial \theta_j} \rangle
 =    \int
   \frac{\partial p(x,\theta)}{\partial \theta_i}\,
   \frac{\partial p(x,\theta)}{\partial \theta_j}\,
   d\lambda(x)
   = h_{ij}(\theta)\ .
\end{eqnarray}
This is different from the Fisher information metric.
The matrix $h(\theta) = (h_{ij}(\theta))$ is called
the direct $L^2$ metric.

Next, we introduce the orthogonal projection between any linear
subspace $V$ of $L^2(\lambda)$
containing the finite dimensional tangent vector space~(\ref{sqrtbas2})
and the tangent vector space~(\ref{sqrtbas2}) itself.
\begin{equation} \label{proL2bis}
\begin{array}{rcl}
  && \Pi_\theta~: L^2(\lambda)\supseteq V \longrightarrow \displaystyle
   \span \{
   \frac{\partial p(\cdot,\theta)}{\partial \theta_1},\cdots,
   \frac{\partial p(\cdot,\theta)}{\partial \theta_m} \} \\ \\
  && \Pi_\theta[v] = \displaystyle \sum_{i=1}^m [ \sum_{j=1}^m
   h^{ij}(\theta)\; \langle v,
   \frac{\partial p(\cdot,\theta)}{\partial \theta_j} \rangle ]\;
   \frac{\partial p(\cdot,\theta)}{\partial \theta_i}\ .
\end{array}
\end{equation}

\begin{example}
{\bf (The Gaussian family and the direct $L^2$ metric in canonical parameters)}. We may consider the $L^2$ metric for the Gaussian family of densities introduced earlier. The $L^2$ metric is
\[ h(\theta) = \frac{1}{8} \frac {\sqrt {2}}{\sqrt {- \theta_2 \pi }} \left[
\begin{array}{cc} 1  &
\frac{\theta_1}{-\theta_2} \\
\frac{\theta_1}{-\theta_2} &
{\frac{3}{4}}\,{\frac{1}{(-\theta_2) }}+ \frac{{\theta_1}^{2}}{{\theta_2}^{2}}
\end{array} \right]
\]
and, as expected, it is different from the Fisher metric seen earlier.
\end{example}

\begin{example}
{\bf (The Gaussian family and the direct $L^2$ metric in expectation parameters)}. We may consider the $L^2$ metric for the Gaussian family in the coordinates $\mu,v$.  The $L^2$ metric is
\[ h(\mu,v) = \frac{1}{8 v \sqrt{v \pi}}  \left[
\begin{array}{cc} 1  &
0 \\
0 & \frac{3}{4v}
\end{array} \right]
\]
and, as expected, it is different from the $\mu,v$ Fisher metric seen earlier, although it is still a diagonal matrix.
\end{example}

%
%
%
%

\section{Exponential families and Mixture families}
Earlier research in \cite{brigo98}, \cite{brigo99}, \cite{brigo99b} and \cite{brigo00} illustrated in detail how the Hellinger distance and the related Fisher information metric are ideal tools when using the projection onto exponential families of densities. This idea was first sketched by Hanzon in \cite{hanzon87}. The above references illustrate this by applying the above framework to the infinite dimensional stochastic PDE describing the optimal solution of the nonlinear filtering problem. This generates an approximate filter that is locally the closest filter in Fisher metric to the optimal one. The use of exponential families allows the correction step in the filtering algorithm to become exact, so that only the prediction step is approximated. Furthermore, and independently from the filtering application, exponential families and the Fisher metric are known to interact well. For example, the Fisher metric is obtained by double differentiation of the normalizing exponent in the exponential family and has a straightforward link with the expectation parameters. See for example \cite{barndorff-nielsen78a}.

The study of the projection filter for exponential families has been carried out in details int he above references, especially \cite{BrigoPhD}, \cite{brigo98} and \cite{brigo99}.

However, besides exponential families, there is another general framework that is powerful in modeling probability densities, and this is the mixture family. Mixture distributions are ubiquitous in statistics and may account for important stylized features such as skewness, multi-modality and fat tails.

We define a mixture family as follows. Suppose we are given $m+1$ fixed squared integrable probability densities in $H_2$, say $\underline{q} = [q_1,q_2,\ldots,q_{m+1}]^T$. Suppose we define the following space of probability densities:

\[ S^M(\underline{q}) = \{ \theta_1 q_1 + \theta_2 q_2 + \cdots + \theta_m q_m + (1-\theta_1 - \cdots - \theta_m) q_{m+1}, \theta_i \ge 0 \ \mbox{for all}\ i, \ \ \theta_1 + \cdots + \theta_m < 1\}\]
 
For convenience, define the transformation

\[ \hat{\theta}(\theta) := [ \theta_1, \theta_2, \ldots, \theta_m, 1 - \theta_1- \theta_2 - \ldots - \theta_m]^T \]
for all $\theta$. We will often write $\hat{\theta}$ instead of $\hat{\theta}(\theta)$ for brevity. With this definition, \[ S^M(\underline{q}) = \{ \hat{\theta}(\theta)^T \uq,\ \   \theta_i \ge 0 \ \mbox{for all}\ i, \ \ \theta_1 + \cdots + \theta_m < 1\}\]
We will occasionally refer to this manifold of densities as to the "Simple Mixture" family. While for exponential families the Hellinger distance and the related Fisher metric are ideal, given also the expression (\ref{Fish1}), for mixture families it is less than ideal. For example, the calculation of the Fisher information matrix $g(\theta)$ becomes cumbersome, and the related projection is quite convoluted. Instead, if we consider the $L^2$ distance and the related structure, the metric $h(\theta)$ and the related projection become very simple. Indeed, one can immediately check from the definition of $h$ that for the mixture family we have
\[ \frac{\partial p(\cdot,\theta)}{\partial \theta_i} = q_i - q_{m+1} \]
and
\[ h_{ij}(\theta) = \int (q_i(x) - q_m(x))(q_j(x) - q_m(x)) d \lambda(x) =:h_{ij}  \]
i.e., the $L^2$ metric (and matrix) does not depend on the specific point $\theta$ of the manifold. The same holds for the tangent space at $p(\cdot,\theta)$, which is given by
\[ L_{p(\cdot,\theta)} S = \mbox{span}\{ q_1-q_{m+1}, q_2-q_{m+1}, \cdots, q_m-q_{m+1}\} \]

Also the $L^2$ projection becomes particularly simple:

\begin{equation} \label{proL2mix}
\begin{array}{rcl}
  && \Pi_\theta~: L^2(\lambda)\supseteq V \longrightarrow \displaystyle
   \mbox{span}\{ q_1-q_{m+1}, q_2-q_{m+1}, \cdots, q_m-q_{m+1}\} \\ \\
  && \Pi_\theta[v] = \displaystyle \sum_{i=1}^m [ \sum_{j=1}^m
   h^{ij}\; \langle v, q_j - q_{m+1}
    \rangle ]\;
   (q_i - q_{m+1})\ .
\end{array}
\end{equation}

It is therefore worthwhile to try and apply the $L^2$ metric and the related structure to the projection of the infinite dimensional filter onto the mixture family above.

\section{The nonlinear filtering problem}
In order to present the key geometric ideas without being overwhelmed by technicalities on stochastic PDEs, we consider the filtering problem with continuous time state and discrete time observations, and in this setup we take a scalar system. We will consider multi-dimensional systems later on, in the case with continuous time observations. 

In this model, the state process is a continuous time stochastic differential equation
\begin{displaymath}
   dX_t = f_t(X_t)\,dt + \sigma_t(X_t)\,dW_t\ ,
\end{displaymath}
but only discrete--time observations are available
\begin{displaymath}
   Z_n = h(X_{t_n}) + V_n\ ,
\end{displaymath}
at times $0 = t_0 < t_1 < \cdots < t_n < \cdots$ regularly sampled,
where $\{V_n\,,\,n\geq 0\}$ is a Gaussian white noise sequence
independent of $\{X_t\,,\,t\geq 0\}$.

The nonlinear filtering problem consists in finding the conditional
density $p_n(x)$ of the state $X_{t_n}$ given the observations
up to time $t_n$,
i.e.\ such that $P[X_{t_n}\in dx\mid {\cal Z}_n] = p_n(x)\,dx$,
where ${\cal Z}_n := \sigma(Z_0,\cdots,Z_n)$.
We define also the prediction conditional density
$p_n^{-}(x)\,dx = P[X_{t_n}\in dx\mid {\cal Z}_{n-1}]$.
The sequence $\{p_n\,,\,n\geq 0\}$ satisfies a recurrent equation,
and the transition from $p_{n-1}$ to $p_n$ is decomposed in two steps,
as explained for example in~\cite{jazwinski70a}.

There is first a prediction step: Between time $t_{n-1}$ and $t_n$, we solve the Fokker--Planck equation
\begin{equation}\label{eq:FPEn}
   \frac{\partial p_t^n}{\partial t} = {\cal L}_t^\ast\, p_t^n\ ,
   \hspace{1cm} p_{t_{n-1}}^n = p_{n-1}\
\end{equation}
where the forward diffusion operator is defined as
\begin{displaymath}
   {\cal L}_t^\ast \phi = -
   \frac{\partial}{\partial x}\, [f_t\, \phi]
   + \half 
   \frac{\partial^2}{\partial x^2}\,
   [\sigma_t^{2}\, \phi]\ ,
\end{displaymath}
while its dual backwards diffusion operator is defined as
\begin{displaymath}
   {\cal L}_t  =  f_t\,
   \frac{\partial}{\partial x} + \half 
   (\sigma_t^2)\, \frac{\partial^2}{\partial x^2 }\ .
\end{displaymath}

The solution at final time $t_n$ defines the prediction
conditional density $p_n^{-} = p_{t_n}^n$.

We have then a second step, the correction step:

At time $t_n$, the newly arrived observation $Z_n$ is combined
with the prediction conditional density $p_n^{-}$
via the Bayes rule
\begin{equation} \label{bayes}
   p_n(x) = c_n\; \Psi_n(x)\; p_n^{-}(x)\ ,
\end{equation}
where $c_n$ is a normalizing constant,
and $\Psi_n(x)$ denotes the likelihood function for the estimation
of $X_{t_n}$ based on the observation $Z_n$ only, i.e.
\begin{equation} \label{likelihood}
   \Psi_n(x) :=
   \exp\open\{ -\half \vert Z_n - h(x) \vert^2 \close\}\ .
\end{equation}

\section{The mixture projection filter (MPF)}
We now introduce the mixture projection filter.

We will now work on the prediction step first, in order to derive the projected version of the Fokker Planck equation, living in the manifold $S^M$.
We adopt the following technique. Take a curve in the mixture family $S^M$,
\[ t \mapsto p(\cdot,\theta(t)) \]
and notice that the left hand side of the Fokker Planck equation for this density would read
\[ \frac{\partial p(\cdot,\theta(t))}{\partial t} = \sum_{i=1}^m \frac{\partial p(\cdot,\theta(t))}{\partial \theta_i}
\frac{d}{dt} \theta_i(t) = \sum_{i=1}^m (q_i-q_{m+1}) \frac{d}{dt} \theta_i(t) \]
and project the right hand side of the Fokker Planck equation as
\[ \Pi_\theta[{\cal L}_t^\ast p(\cdot,\theta)] = \sum_{i=1}^m [ \sum_{j=1}^m
   h^{ij}\; \langle {\cal L}_t^\ast p(\cdot,\theta), q_j - q_{m+1}
    \rangle ]\;
   (q_i - q_{m+1})\ =\]
\[ = \sum_{i=1}^m [ \sum_{j=1}^m
   h^{ij}\; \langle p(\cdot,\theta), {\cal L}_t( q_j - q_{m+1})
    \rangle ]\;
   (q_i - q_{m+1})\]
where we used integration by parts in the last step.
Now equating the two sides we obtain
\[ \sum_{i=1}^m (q_i-q_{m+1}) \frac{d}{dt} \theta_i(t) = \sum_{i=1}^m [ \sum_{j=1}^m
   h^{ij}\; \langle p(\cdot,\theta), {\cal L}_t( q_j - q_{m+1})
    \rangle ]\;
   (q_i - q_{m+1})\]
which yields the ordinary differential equation for the parameters $\theta$ of the projected density:
\[ \frac{d}{dt} \theta_i(t) = \sum_{j=1}^m
   h^{ij}\; \langle p(\cdot,\theta), {\cal L}_t( q_j - q_{m+1})
    \rangle \]
Now, by taking into account the structure of $p(\cdot,\theta)$ and the fact that such densities are linear in $\theta$, we see that the above equation is a linear differential equation:

\[ \frac{d}{dt} \theta_i(t) = \sum_{j=1}^m
   h^{ij}\; \left[ \sum_{k=1}^m \theta_k \langle q_k, {\cal L}_t( q_j - q_{m+1})
    \rangle + (1-\theta_1 - \cdots - \theta_m )\langle q_{m+1},
{\cal L}_t( q_j - q_{m+1})\rangle \right] . \]

If we define, for two vector functions $f$ and $g$, the matrix $\langle f, g \rangle$ and the vector
${\cal L}_t f$
as
\[ (\langle f, g \rangle)_{i,j} := \langle f_i, g_j \rangle, \ \ ({\cal L}_t f)_i := {\cal L}_t (f_i) \]
then we can write the above ODE in compact form as
\begin{equation}\label{eq:ProFPE} 
\frac{d}{dt} \underline{\theta}(t) = h^{-1} \langle\ {\cal L}_t( \uq_{1:m} - 1_m q_{m+1}) , \ \uq\ \rangle \ \hat{\theta}(\underline{\theta}(t)) 
\end{equation}
where $\uq_{1:m}$ is the vector with the first $m$ components of $\uq$, and $1_m$ is a m-dimensional (column) vector of ones.

In \cite{brigo98} and \cite{brigo99} it is shown that, by carefully choosing the exponential family, the Fisher metric exponential projection filter makes the correction step exact. In the mixture framework under the $L^2$ metric we are using now, this is harder to achieve unless we are willing to redefine the manifold at every correction step. Let us therefore focus on the correction step first. Suppose we are in $[t_{n-1},t_n)$ and we obtained a prediction for the density up to $t_n^-$, whose parameter we call $\underline{\theta}_n^-$. At $t_n$ a new observation $Z_n$ arrives and we update the density. Substituting the prediction $p(\cdot,\underline{\theta}_n^{-})$ into formula~(\ref{bayes}), we observe that the resulting density leaves the original mixture family $S^M(\underline{q})$.
The updated density at $t_n$ is
\[ c_n \Psi_n(x) p(x,\theta_n^{-}) = c_n \Psi_n(x)\ \hat{\underline{\theta}}^T \underline{q}    \]
and is outside $S^M(\underline{q})$. However, we may keep the update step exact by re-defining the basis functions $q$ as follows.

Suppose that we change basis functions at every discrete date observation step. The first basis function vector is
$\uq^0$, then at update time $t_1$ we will select a new vector of basis functions $\uq^1$, and so on. At every point in time we keep the vector $m+1$ dimensional. Suppose the basis functions in $[t_{n-1}, t_{n})$ are $\uq^{n-1}$.
We run the prediction step up to $t_{n}^-$, getting $\underline{\theta}_n^-$.  At time $t_n$, we define the new basis functions as
\[ q^{n}_i(x) := c_{i,n} \Psi_n(x) q^{n-1}_i(x) \ \ \mbox{for all} \ i=1,\ldots,m+1 \]
and where $c_{i,n}$ is the normalizing constant for the density on the right hand side.
Every $q^{n}_i$ is a normalized densities and we can define a mixture of such densities as the new space.
In this case, the correction step amounts to set, at $t_n$:

\paragraph{Correction Step:}
\[\mbox{At} \ t_n: \ \underline{\theta}_n = \underline{\theta}^n_{t_n}, \ \ \mbox{and the new manifold is} \ S^M(\uq^{n})  \]

We may now focus on the prediction step.

Before doing so, it is important to notice that the $L^2$ metric changes as well when we change the manifold, so that it is safe to index as follows:

\[ h^n_{ij} = \int (q^n_i(x) - q^n_m(x))(q^n_j(x) - q^n_m(x)) d \lambda(x)  \]

\paragraph{Prediction step}
Between time $t_{n-1}$ and $t_n$, we solve the ODE's
\[ \frac{d}{dt} \underline{\theta}^n(t) = (h^{n-1})^{-1} \langle {\cal L}_t( \uq^{n-1}_{1:m} - 1_m q^{n-1}_{m+1}), \ \uq^{n-1} \rangle \ \hat{\theta}(\underline{\theta}^n(t)), \ \  \underline{\theta}_{t_{n-1}}^n := \underline{\theta}_{n-1}\ . \]
The solution at final time $t_n$ defines the prediction
parameters $\underline{\theta}_n^{-} = \theta_{t_n}^n$.

\section{Relationship with Galerkin methods}
Consider the prediction step in the above section. 
This is the same step we would have obtained through Galerkin methods, see for example \cite{beard}. 

In \cite{beard}, the Galerkin method is applied to the filtering problem with continuous time observations. We will address the continuous time observations  setup in the next section. 

Here we keep discrete time observations and we show the Galerkin approximation on the prediction step.

The Galerkin approximation is obtained by approximating the exact solution of the Fokker--Planck equation (\ref{eq:FPEn}) with a function of the form
\begin{equation}\label{eq:galerkin:approx} \tilde{p}_t(x) := \sum_{i=1}^{m+1} c_i(t) \phi_i(x) , 
\end{equation}
see for example \cite{beard} for more details. The method works by replacing the exact solution of the Fokker--Planck equation with the solution of the equations
\[   \langle - \frac{\partial \tilde{p}_t}{\partial t} + {\cal L}_t^\ast\, \tilde{p}_t\, , \xi \rangle = 0 \]
for a suitable family of smooth $L^2$ test functions $\xi$. 
By using the approximation (\ref{eq:galerkin:approx}) in this last expression, and by taking $\xi = \phi_j$ for $j=1,\ldots,m+1$, and finally by setting
\[ c_i(t) = \theta_i(t)\ \mbox{and} \ \phi_i(x) = q_i(x) - q_{m+1}(x) \ \mbox{for} \ i=1,\ldots,m,\ \  c_{m+1}(t)=1,  \ \phi_{m+1}(x) = q_{m+1}(x)\]
we can see that the method provides exactly Equation (\ref{eq:ProFPE}). Therefore, for simple mixture families the $L^2$ projection filter prediction step will coincide with the Galerkin method based prediction step.  

However, this holds only for the case where the manifold $S$ on which we project is the simple mixture family. More complex families, such as the ones we will hint at in the continuous observation case, will not allow for a Galerkin-based filter and only the $L^2$ projection filter can be defined there. Furthermore,  even under the simple mixture family, in the continuous observations case there is a further fundamental difference. Our $L^2$ projection filter in the continuos time observations case will be different from the Galerkin projection filter in \cite{beard}, because we use Stratonovich calculus to project the Kushner-Stratonovich equation in $L^2$ metric. In \cite{beard} the Ito version of the Kushner-Stratonovich Equation is used instead, but since Ito calculus does not work on manifolds, due to the second order term moving the dynamics out of the tangent space (see for example \cite{brigo99b}), we use the Stratonovich version instead. The Ito-based and Stratonovich based Galerkin projection filters will therefore differ for simple mixture families, and again, only the second one can be defined for manifolds of densities beyond the simplest mixture family. A particularly important manifold for which only the $L^2$ based filter can be defined is a manifold that makes the correction step exact also in continuous time.  For such a family one can define a rigorous measure of the filtering error in $L^2$ norm, which is impossible to obtain with the standard Galerkin method. This will be made explicit in future work.

\section{The Filtering Problem with continuous-time observations }

In the above part of the paper we decided to take discrete time observations in order to limit technicalities. In this section we consider a continuos time framework both for the observations $Y$ and for the signal $X$, and we allow both to be multi-dimensional processes.

\begin{equation} \label{Lanc1-1}
\begin{array}{rcl}
   dX_t &=& f_t(X_t)\,dt + \sigma_t(X_t)\,dW_t, \ \ X_0, \\ \\
   dY_t &=& b_t(X_t)\,dt + dV_t, \ \ Y_0 = 0\ .
\end{array}
\end{equation}
These equations are It\^o stochastic differential equations (SDE's).
In the continuous observations case we shall use both It\^o SDE's (for example for the
signal $X$) and Stratonovich (Str) SDE's
(when dealing with manifolds and projections). The Str form
will be denoted by the presence of the symbol `$\circ$' in between
the diffusion coefficient and the Brownian motion of a SDE.
The use of Str SDE's is necessary in order to be able
to deal with stochastic calculus on manifolds, since in general
one does not know how to interpret the second order terms arising in
It\^o's calculus in terms of manifold structures.
The interested reader is referred to~\cite{elworthy82a}.

In~(\ref{Lanc1-1}), the unobserved state process $\{X_t\,,\,t\geq 0\}$
and the observation process $\{Y_t\,,\,t\geq 0\}$ are taking
values in ${\bf R}^n$ and ${\bf R}^d$ respectively,
the noise processes $\{W_t\,,\,t\geq 0\}$ and $\{V_t\,,\,t\geq 0\}$
are two Brownian motions, taking values in ${\bf R}^p$ and ${\bf R}^d$
respectively, with covariance matrices $Q_t$ and $R_t$ respectively.
We assume that $R_t$ is invertible for all $t\geq 0$,
which implies that, without loss of generality,
we can assume that $R_t = I$ for all $t\geq 0$.
Finally, the initial state $X_0$ and the noise processes
$\{W_t\,,\,t\geq 0\}$ and $\{V_t\,,\,t\geq 0\}$ are assumed to be
independent.
We assume that the initial state $X_0$ has a density $p_0$
w.r.t.\ the Lebesgue measure $\lambda$ on ${\bf R}^n$,
and has finite moments of any order,
and we make the following assumptions
on the coefficients $f_t$, $a_t := \sigma_t\, Q_t\, \sigma_t^T$,
and $b_t$ of the system~(\ref{Lanc1-1})
\begin{itemize}
   \item[(A)] Local Lipschitz continuity~:
for all $R > 0$, there exists $K_R > 0$ such that
\begin{displaymath}
   \vert f_t(x) - f_t(x') \vert \leq K_R\, \vert x-x' \vert
   \hspace{1cm}\mbox{and}\hspace{1cm}
   \Vert a_t(x) - a_t(x') \Vert \leq K_R\, \vert x-x' \vert\ ,
\end{displaymath}
for all $t\geq 0$, and for all $x,x'\in B_R$, the ball of radius $R$.

   \item[(B)] Non--explosion~:
there exists $K > 0$ such that
\begin{displaymath}
   x^T f_t(x) \leq K\, (1+\vert x \vert^2)
   \hspace{1cm}\mbox{and}\hspace{1cm}
   \mbox{trace}\; a_t(x) \leq K\, (1+\vert x \vert^2)\ ,
\end{displaymath}
for all $t\geq 0$, and for all $x\in {\bf R}^n$.

   \item[(C)] Polynomial growth~:
there exist $K > 0$ and $r \geq 0$ such that
\begin{displaymath}
   \vert b_t(x) \vert \leq K\, (1+\vert x \vert^r)\ ,
\end{displaymath}
for all $t\geq 0$, and for all $x\in {\bf R}^n$.
\end{itemize}

Under assumptions~(A) and~(B), there exists
a unique solution $\{X_t\,,\,t\geq 0\}$ to the state equation,
see for example~\cite{khasminskii},
and $X_t$ has finite moments of any order.
Under the additional assumption~(C) the following {\em finite energy}
condition holds
\begin{displaymath}
   E \int_0^T \vert b_t(X_t) \vert^2\,dt < \infty\ ,
   \hspace{1cm}\mbox{for all $T\geq 0$}.
\end{displaymath}

The nonlinear filtering problem consists in finding
the conditional probability distribution $\pi_t$ of the state $X_t$
given the observations up to time $t$,
i.e.\ $\pi_t(dx) := P[X_t\in dx\mid {\cal Y}_t]$,
where ${\cal Y}_t:=\sigma(Y_s\,,\,0\leq s\leq t)$.
Since the finite energy condition holds,
it follows from Fujisaki, Kallianpur and Kunita~\cite{fujisaki72a}
that $\{\pi_t\,,\,t\geq 0\}$ satisfies the Kushner--Stratonovich equation,
i.e.\ for any smooth and compactly supported test function $\phi$
defined on ${\bf R}^n$
\begin{equation} \label{FKK}
   \pi_t(\phi) = \pi_0(\phi)
   + \int_0^t \pi_s({\cal L}_s \phi)\,ds
   + \sum_{k=1}^d \int_0^t
   [\pi_s(b_s^k\, \phi) - \pi_s(b_s^k)\, \pi_s(\phi)]\,
   [dY_s^k-\pi_s(b_s^k)\,ds]\ ,
\end{equation}
where for all $t \geq 0$,
the backward diffusion operator ${\cal L}_t$ is defined by
\begin{displaymath}
   {\cal L}_t = \sum_{i=1}^n f_t^i\,
   \frac{\partial}{\partial x_i} + \half \sum_{i,j=1}^n
   a_t^{ij}\, \frac{\partial^2}{\partial x_i \partial x_j}\ .
\end{displaymath}
The Str form of equation~(\ref{FKK}) is obtained, after straightforward
computations, as~:
\begin{equation} \label{MFS-FKK}
\begin{array}{rcl}
   \pi_t(\phi) &=& \displaystyle \pi_0(\phi)
   + \int_0^t \pi_s({\cal L}_s\, \phi)\,ds
   - \half\, \int_0^t [\pi_s(\vert b_s \vert^2\, \phi)
   - \pi_s(\vert b_s \vert^2)\, \pi_s(\phi)]\,ds \\ \\
   && \displaystyle \qquad \mbox{}
   + \sum_{k=1}^d \int_0^t
   [\pi_s(b_s^k\, \phi) - \pi_s(b_s^k)\, \pi_s(\phi)]\circ dY_s^k\ .
\end{array}
\end{equation}
From now on we proceed formally, and we assume that for all $t \geq 0$,
the probability distribution $\pi_t$ has a density $p_t$
w.r.t.\ the Lebesgue measure on ${\bf R}^n$.
Then $\{p_t\,,\,t\geq 0\}$ satisfies the It\^o--type stochastic
partial differential equation (SPDE)
\begin{equation} \label{Kushner}
   dp_t = {\cal L}_t^\ast\, p_t\,dt
   + \sum_{k=1}^d p_t\, [b_t^k-E_{p_t}\{b_t^k\}]\,
   [dY_t^k-E_{p_t}\{b_t^k\}\,dt]
\end{equation}
in a suitable functional space,
where $E_{p_t}\{\cdot\}$ denotes the expectation w.r.t.\
the probability density $p_t$, i.e.\ the conditional expectation
given the observations up to time $t$,
and where for all $t \geq 0$,
the forward diffusion operator ${\cal L}_t^\ast$ is defined by
\begin{displaymath}
   {\cal L}_t^\ast \phi = -\sum_{i=1}^n
   \frac{\partial}{\partial x_i}\, [f_t^i\, \phi]
   + \half \sum_{i,j=1}^n
   \frac{\partial^2}{\partial x_i \partial x_j}\,
   [a_t^{ij}\, \phi]\ ,
\end{displaymath}
for any test function $\phi$ defined on ${\bf R}^n$.
%
%
The corresponding Str form of the SPDE~(\ref{Kushner}) is~:
\begin{displaymath}
   dp_t = {\cal L}_t^\ast\, p_t\,dt
   - \half\, p_t\, [\vert b_t \vert^2 - E_{p_t}\{\vert b_t \vert^2\}] \,dt
   + \sum_{k=1}^d p_t\, [b_t^k-E_{p_t}\{b_t^k\}] \circ dY_t^k\ .
\end{displaymath}
In order to simplify notation, we introduce the following
definitions~:
\begin{equation} \label{coeff}
\begin{array}{rcl}
 \gamma_t^0(p) &:=&
   \half\, [\vert b_t \vert^2 - E_p\{\vert b_t \vert^2\}]\ p, \\ \\
   \gamma_t^k(p) &:=& [b_t^k - E_p\{b_t^k\}] p \ ,
\end{array}
\end{equation}
for $k = 1,\cdots,d$.
The Str form of the Kushner--Stratonovich equation reads now
\begin{equation}\label{KSE:str}
   dp_t = {\cal L}_t^\ast\, p_t\,dt
   -  \gamma_t^0(p_t)\,dt
   + \sum_{k=1}^d  \gamma_t^k(p_t) \circ dY_t^k\ .
\end{equation}

This equation can be projected according to the L2 direct metric we introduced above, similarly to how we projected the Fokker Planck equation for the prediction step in the discrete time observation case. There the projection transformed a PDE into a ODE, whereas in our current case the projection will transform a SPDE into a SDE.

Take again a curve in the mixture family $S^M$,
\[ t \mapsto p(\cdot,\theta(t)) \]
and notice that the left hand side of the Kushner-Strantonovich SPDE for this density would read
\[ d_t p(\cdot,\theta(t)) = \sum_{i=1}^m \frac{\partial p(\cdot,\theta(t))}{\partial \theta_i}
d_t \theta_i(t) = \sum_{i=1}^m (q_i-q_{m+1})\ d \theta_i(t) \]
and project the right hand side terms of the Kushner-Strantonovich SPDE (\ref{KSE:str}) as
\[ \Pi_\theta[{\cal L}_t^\ast p(\cdot,\theta)]  = \sum_{i=1}^m [ \sum_{j=1}^m
   h^{ij}\; \langle p(\cdot,\theta), {\cal L}_t( q_j - q_{m+1})
    \rangle ]\;
   (q_i - q_{m+1}),\]
\[ \Pi_\theta[\gamma^k_t(p(\cdot,\theta)) ] = \sum_{i=1}^m [ \sum_{j=1}^m
   h^{ij}\; \langle \gamma^k_t(p(\cdot,\theta)), q_j - q_{m+1}
    \rangle ]\;
   (q_i - q_{m+1})\ \]
Now equating the two sides we obtain
\[ \sum_{i=1}^m (q_i-q_{m+1}) {d_t} \theta_i(t) = \sum_{i=1}^m \bigg{[} \sum_{j=1}^m
   h^{ij}\;\bigg{\{} \langle p(\cdot,\theta), {\cal L}_t( q_j - q_{m+1})
    \rangle dt - \langle \gamma^0_t(p(\cdot,\theta)), q_j - q_{m+1}
    \rangle dt  \]\[+ \sum_{k=1}^d  \langle \gamma^k_t(p(\cdot,\theta)), q_j - q_{m+1}
    \rangle  \circ dY_t^k \bigg{\}}  \bigg{]}\;
   (q_i - q_{m+1})\]
which yields the stochastic differential equation for the parameters $\theta$ of the projected density:
\[ d_t \theta_i(t) =  \sum_{j=1}^m
   h^{ij}\;\bigg{\{} \langle p(\cdot,\theta), {\cal L}_t( q_j - q_{m+1})
    \rangle dt - \langle \gamma^0_t(p(\cdot,\theta)), q_j - q_{m+1}
    \rangle dt  \]\[+ \sum_{k=1}^d  \langle \gamma^k_t(p(\cdot,\theta)), q_j - q_{m+1}
    \rangle  \circ dY_t^k \bigg{\}}  \;
   \]
Similarly to what we did for the discrete time observations case, we can write this SDE in more compact form as
\[ {d_t} \underline{\theta}(t) = h^{-1} \langle {\cal L}_t( \uq_{1:m} - 1_m q_{m+1}), \ \uq \rangle \ \hat{\theta}(\underline{\theta}(t)) dt - h^{-1} \langle \gamma^0_t(p(\cdot,\theta)), \uq_{1:m} - 1_m q_{m+1}) \rangle dt+\]\begin{equation}\label{eq:PFcto}
+ h^{-1} \sum_{k=1}^d  \langle \gamma^k_t(p(\cdot,\theta)), \uq_{1:m} - 1_m q_{m+1}  \rangle  \circ dY_t^k  \end{equation}

Notice that now only the prediction $dt$ part is linear. More generally, by inspection one can see that the equation is quadratic. One can define a projection residual in $L^2$ norm, measuring the local projection error of the filter. This residual can be made rigorous under a specific mixture family incorporating a pseudo likelihood ratio update factor into each mixture family member function $q$. This, and a numerical investigation on the effectiveness of the filter for some standard systems is under investigation in \cite{brigocapponipf}.

\section{Conclusion and Further Research}
We introduced a projection method for the space of probability distributions based on the differential geometric approach to statistics. This method makes use of a direct $L^2$ metric as opposed to the usual Hellinger distance and the related Fisher Information metric. We applied this apparatus to the nonlinear filtering problem. Past projection filters concentrated on the Fisher metric and the exponential families that made the filter correction step exact. Instead, in this work we introduce the mixture projection filter, namely the projection filter based on the direct $L^2$ metric and based on a manifold given by a mixture of pre-assigned densities. We derived the filter equations for the discrete time observation case first. We showed how an update on the manifold functions, even when keeping the same dimension, can make the correction step exact. A key result is that the prediction step is a simple linear ordinary differential equation.

We then derived the continuous time observations filter by projecting the Kushner Stratonovich stochastic PDE in Stratonovich form, and obtained a SDE whose drift is linear but with additional quadratic terms both in the drift and in the diffusion part.

We finally remarked that the exponential projection filter had a clear relationship with the assumed density filters, as documented in \cite{brigo99}. The mixture projection filter introduced here has a clear relationship with earlier Galerkin-based approaches when applied to simple mixture families, although even for such families there are important differences in the continuous time observations case. In future work \cite{brigocapponipf} we will also implement the mixture projection filter equations numerically and will investigate the choice of the specific mixture family and the projection error.

\end{document}